\documentclass[12pt, reqno]{amsart}

\usepackage[hmargin=0.8in,twosideshift=0in,height=8.6in]{geometry}
\usepackage{amssymb,amsthm}
\usepackage{delarray,verbatim}
\usepackage{natbib}

\usepackage{ifpdf}
\ifpdf
\usepackage[pdftex]{graphicx}
 \else
\usepackage[dvips]{graphicx}
 \fi

\usepackage{ifpdf}
\usepackage{color}
\definecolor{webgreen}{rgb}{0,.5,0}
\definecolor{webbrown}{rgb}{.8,0,0}
\definecolor{emphcolor}{rgb}{0.95,0.95,0.95}

\usepackage{hyperref}
\hypersetup{%
          colorlinks=true,
          linkcolor=webbrown,
          filecolor=webbrown,
          citecolor=webgreen,
          breaklinks=true}
\ifpdf \hypersetup{pdftex,
            pdfstartview=FitH, 
            bookmarksopen=true,
            bookmarksnumbered=true
} \else \hypersetup{dvips} \fi

\usepackage{multicol}
\usepackage{multirow}

\newtheorem{theorem}{Theorem}

\newcommand{\E}{\mathbb{E}}

\newcommand{\R}{\mathbb{R}}
\newcommand{\N}{\mathbb{N}}

\title[]{On the uniqueness of classical solutions of Cauchy problems}\thanks{This research is supported in part by the National Science Foundation under grant number DMS-0906257.}

\author[]{Erhan Bayraktar}
\address[Erhan Bayraktar]{Department of Mathematics, University of Michigan, 530 Church Street, Ann Arbor, MI 48104, USA}
\email{erhan@umich.edu}
\author[]{Hao Xing }
\address[Hao Xing]{Department of Mathematics and Statistics, Boston University, Boston, MA 02215}
\email{haoxing@umich.edu}

\date{}
\begin{document}

\begin{abstract}
Given that the terminal condition is of at most linear growth, it is well known that a Cauchy problem admits a unique classical solution when the coefficient multiplying the second derivative is also a function of at most linear growth. In this note, we give a condition on the volatility that is necessary and sufficient for a Cauchy problem to admit a unique solution 

\end{abstract}

\maketitle

\section{Main Result}
Given a terminal-boundary data $g:\R_+ \mapsto \R_+$ with $g(x) \leq C(1+x)$ for some constant $C>0$, we consider the following Cauchy problem
\begin{equation}\label{eq:Cauchy problem}
\begin{split}
 & u_t + \frac12 \sigma^2(x) \, u_{xx} = 0, \quad (x,t) \in (0,\infty) \times [0,T), \\
 & u(0, t) = g(0), \quad t\leq T,\\
 & u(x, T) = g(x),
\end{split}
\end{equation}
where $\sigma \neq 0$ on $(0,\infty)$, $\sigma^{-2} \in L_{loc}^1(0,\infty)$ (i.e., $\int_{a}^{b}\sigma^{-2}(x)dx<\infty$ for any $[a,b] \subset (0,\infty)$), and $\sigma=0$ on $(-\infty,0]$. 

A solution $u: \R_+ \times [0, T] \mapsto \R$ of \eqref{eq:Cauchy problem} is said to be a \emph{classical solution} if  $u\in C^{2,1} ((0,\infty) \times [0,T))$. A function $f: \R_+ \times [0,T] \mapsto \R$ is said to be of \emph{at most linear growth} if there exists a constant $C>0$ such that $|f(x,t)|\leq C(1+x)$ for any $(x,t) \in \R_+ \times [0,T]$. 

A well-known sufficient condition for \eqref{eq:Cauchy problem} to have a unique classical solution among the functions with at most linear growth is that $\sigma$ itself is of at most linear growth; see e.g. Chapter 6 of \cite{MR0494490} and Theorem 7.6 on page 366 of \cite{MR1121940}. On the other hand, consider the SDE
 \begin{equation}\label{eq:SDE}
  d X^{t,x}_s = \sigma(X^{t,x}_s) \, dW_s, \quad X^{t,x}_t = x>0.
 \end{equation}
The assumptions on $\sigma$ we made below \eqref{eq:Cauchy problem} ensure that \eqref{eq:SDE} has a unique weak solution which is absorbed at zero. (See \cite{MR798317}.) The solution $X^{t,x}$ is clearly a local martingale. Delbaen and Shirakawa shows in \cite{citeulike:64915} that $X^{t,x}$ is a martingale if and only if 
\begin{equation}\label{condition}
  \int_1^{\infty} \frac{x}{\sigma^2(x)}  \, dx = \infty.
\end{equation}
Also see \cite{CCU2009}.
 
 Below, in Theorem~\ref{thm:main2}, we prove that \eqref{condition} is also necessary and sufficient for the existence of a unique classical solution of \eqref{eq:Cauchy problem}. First, in the next theorem, we show that \eqref{condition}, which is weaker than the linear growth condition on $\sigma$, is a sufficient condition for the uniqueness.

\begin{theorem} \label{thm:main1}
 The Cauchy problem \eqref{eq:Cauchy problem} has a unique classical solution (if any) in the class of functions with at most linear growth if
 \eqref{condition} is satisfied.
\end{theorem}

\begin{proof}
It suffices to show that $u\equiv 0$ is the unique solution of the Cauchy problem in \eqref{eq:Cauchy problem} with $g\equiv 0$. 
Let us define a sequence of stopping times $\tau_n \triangleq \inf\{s \geq t : X^{t,x}_s \geq n\; \text{or}\; X^{t,x}_s \leq 1/n\} \wedge T$ for each $n\in \N_+$ and $\tau_0 \triangleq \{s\geq t : X_s^{t,x} = 0\} \wedge T$. Then as in the proof of Theorem 1.6 in \cite{citeulike:64915} we can show that the function defined by
\[
\Psi(x)=
\begin{cases}
x, & x \leq 1;
\\ x+\int_1^{x}\frac{u}{\sigma^2(u)}(x-u)du, & x \geq 1.
\end{cases}
\]
satisfies $\E[\Psi(X^{t,x}_{\tau_n})] \leq \Psi(x) + xT /2$. Since $\Psi$ is convex, \eqref{condition} implies that $\lim_{x \to \infty}\Psi(x)/x=\infty$. Then the criterion of de la Vall\'{e}e Poussin implies that $\{X^{t,x}_{\tau_n} : n \in \N\}$ is a uniformly integrable family.

 Suppose $\tilde{u}$ is another classical solution of at most linear growth. Applying the It\^{o}'s lemma, we obtain
\[
 \tilde{u}(X^{t,x}_{\tau_n}, \tau_n) = \tilde{u}(x, t) + \int_t^{\tau_n} \left[ \tilde{u}_s (X^{t,x}_s, s) + \frac12 \sigma^2(X^{t,x}_s) \tilde{u}_{xx} (X^{t,x}_s, s)\right] \, ds + \int_t^{\tau_n} \tilde{u}_x (X^{t,x}_s, s) \sigma(X^{t,x}_s) dW_s.
\]
Thanks to our choice of $\tau_n$, the expectation of the stochastic integral is zero. Therefore, taking the expectation of both sides of the above identity, we get $\tilde{u}(x,t) = \E\left[\tilde{u}(X^{t,x}_{\tau_n}, \tau_n)\right]$ for each $n\in \N$.

On the other hand, since $\tilde{u}$ is of at most linear growth, there exists a constant $C$ such that $|\tilde{u}(x,t) | \leq C(1+x)$. Therefore $\{\tilde{u}(X^{t,x}_{\tau_n}, \tau_n) : n\in \N_+\}$ is a uniformly integrable family. This is because it is bounded from above by the uniformly integrable family $\{C (1+X^{t,x}_{\tau_n}):n\in \N_+\}$. As a result,
\[
\begin{split}
 \tilde{u}(x,t) &=\lim_{n\rightarrow \infty} \E\left[\tilde{u}(X^{t,x}_{\tau_n}, \tau_n)\right] = \E\left[\lim_{n\rightarrow \infty} \tilde{u}(X^{t,x}_{\tau_n}, \tau_n)\right] = \E[\tilde{u} (X^{t,x}_{\tau_0}, \tau_0)] \\
 & =\E\left[g(X_T^{t,x}) 1_{\{\tau_0 =T\}}\right] + \E\left[\tilde{u}\left(X_{\tau_0}^{t,x}, \tau_0\right) 1_{\{\tau_0 <T\}}\right] \\
 & = \E\left[g(X_T^{t,x}) 1_{\{\tau_0 =T\}}\right] + \E\left[g(0) 1_{\{\tau_0 <T\}}\right] =0.
 \end{split}
\]
Here the third equality holds since $X^{t,x}$ does not explode (i.e., $\inf\{s\geq t : X_s^{t,x} = +\infty\} = \infty$, see Problem 5.3 in page 332 of \cite{Karatzas-Shreve}) and one before the last equality follows since $X^{t,x}_{\tau_0}=0$ on the set $\{\tau_0 < T\}$.
\end{proof}

\begin{theorem}\label{thm:main2}
If we further assume that $\sigma : \R_+ \mapsto \R_+$ is locally H\"{o}lder continuous with exponent $1/2$ and $g$ is of linear growth, then the Cauchy problem in \eqref{eq:Cauchy problem} has a unique classical solution if and only if \eqref{condition} is satisfied.
\end{theorem}

\begin{proof}

First let us prove the existence of a solution.
Let  $u(x,t) \triangleq \E \, g(X^{t,x}_T)$ (the value of a call-type European option). Thanks to the H\"{o}lder continuity of $\sigma$, it follows from Theorem 3.2 in \cite{Ekstrom-Tysk} that $u$ is a classical solution of \eqref{eq:Cauchy problem}. Moreover, it is of at most linear growth due to the assumption that $g$ is of at most linear growth.

\emph{Proof of sufficiency.} This follows from Theorem~\ref{thm:main1}.

 \emph{Proof of necessity.}  If \eqref{condition} is violated, then $X$ is a strict local martingale (see \cite{citeulike:64915} and \cite{CCU2009}). Using Theorem 3.2 in \cite{Ekstrom-Tysk}, it can be seen that $u^*(x,t) \triangleq x - \E[X^{x,t}_T] >0$ is a classical solution of \eqref{eq:Cauchy problem} with zero boundary and terminal conditions. (Note that the H\"{o}lder continuity assumption on $\sigma$ is used in this step as well.) This function clearly has at most linear growth. Therefore $u + \lambda u^*$, for any $\lambda \in \R$, is also a classical solution of \eqref{eq:Cauchy problem} which is of at most linear growth.
\end{proof}

A related result is given by Theorem 4.3 of  \cite{Ekstrom-Tysk} on put-type European options: When $g$ is of strictly sublinear growth (i.e., $\lim_{x\rightarrow \infty} g(x)/x =0$) then \eqref{eq:Cauchy problem} has a unique solution among the functions with strictly sublinear growth (without assuming \eqref{condition}).

Our result in Theorem~\ref{thm:main2} complements Theorem 3.2 of \cite{Ekstrom-Tysk}, which shows that the call-type European option price is a classical solution of \eqref{eq:Cauchy problem} of at most linear growth. We prove that \eqref{condition} is necessary and sufficient to guarantee that the European option price is the only classical solution (of at most linear growth) to this Cauchy problem. \cite{Cox-Hobson} and \cite{Heston-Loewenstein-Willard} had already observed that the Cauchy problem corresponding to European call options have multiple solutions. (Also see  \cite{DFern-Kar09} and \cite{bkx09}, which consider super hedging prices of call-type options when there are no equivalent local martingale measures.) However, a necessary and sufficient condition under which there is uniqueness/nonuniqueness remained unknown.

\bibliographystyle{plain}
\bibliography{biblio}

\end{document}